# Numerical Approximations of a Class of Nonlinear Second-Order Boundary Value Problems using Galerkin-Compact Finite Difference Method

Shovan Sourav Datta Pranta, and Md. Shafiqul Islam*

Department of Applied Mathematics, University of Dhaka, Bangladesh. (e-mail: shovan275@gmail.com)
*Corresponding author and e-mail: mdshafiqul@du.ac.bd

**ABSTRACT**

In this study, we examine numerical approximations for 2nd-order linear-nonlinear differential equations with diverse boundary conditions, followed by the residual corrections of the first approximations. We first obtain numerical results using the Galerkin weighted residual approach with Bernstein polynomials. The generation of residuals is brought on by the fact that our first approximation is computed using numerical methods. To minimize these residuals, we use the compact finite difference scheme of 4th-order convergence to solve the error differential equations in accordance with the error boundary conditions. We also introduce the formulation of the compact finite difference method of fourth-order convergence for the nonlinear BVPs. The improved approximations are produced by adding the error values derived from the approximations of the error differential equation to the weighted residual values. Numerical results are compared to the exact solutions and to the solutions available in the published literature to validate the proposed scheme, and high accuracy is achieved in all cases.

## I. INTRODUCTION

The capacity to provide the nature of the solution to any physical occurrence, even when analytical answers are not possible, is a key advantage of the numerical approach. In addition, a numerical method requires simply the evaluation of standard functions and the four operations of addition, subtraction, multiplication, and division. There has been a significant increase in interest in the study of boundary value problems for second- and higher-order differential equations because both linear and nonlinear differential equations have the ability to replicate a wide variety of natural processes. They are also implemented in a variety of scientific and engineering applications.

In this study, we consider the following generic form of a second-order linear-nonlinear boundary value problem:

$$\frac{d^2 f}{dx^2} = g\left(x, f, \frac{df}{dx}\right) \qquad a \leq x \leq b \tag{1}$$

in accordance with the boundary conditions,

$$\left. \begin{array}{l} a_1 f(a) + b_1 \dfrac{df}{dx}\bigg|_{x=a} = \lambda_1 \\[2mm] a_2 f(b) + b_2 \dfrac{df}{dx}\bigg|_{x=b} = \lambda_2 \end{array} \right\} \tag{2}$$

where $a_1, a_2, b_1, b_2, \lambda_1, \lambda_2$ are some constants. If
1. $b_1 = b_2 = 0$ and $a_1 \neq 0, a_2 \neq 0$, (2) is referred to as Dirichlet boundary conditions.
2. $a_1 = a_2 = 0$ and $b_1 \neq 0, b_2 \neq 0$, (2) is referred to as Neumann boundary conditions.
3. If all of the parameters $a_1, a_2, b_1, b_2$ are not equal to 0, the Robin boundary condition is established.

In the field of numerical analysis, there are a lot of different ways to deal with boundary value problems (BVPs). In their respective books, Bender *et al.* [1] and Collatz [2] provided a full theory and some basic numerical treatments of boundary value problems. Kanth and Reddy [3] employed the cubic spline technique for solving two-point boundary value problems numerically. Then using the same technique, Kanth and Bhattacharya [4] analyzed and numerically resolved a class of nonlinear boundary value problems appearing in physiology. Sibana *et al.* [5] used a new spectral-homotopy model to approximate the solution of second-order nonlinear BVPs. A class of second-order differential equations was solved numerically by Islam and Shirin [6] using a weighted residual technique called the Galerkin method via Bernoulli polynomials. In certain instances, however, a significant number of Bernoulli polynomials are necessary to achieve the needed precision, necessitating a substantial amount of computational time. Burden and Faires [7] explored certain numerical techniques for handling boundary value problems in their book on numerical analysis, such as the shooting method and the finite difference method. While, using the parametric difference approaches, Pandey [8] solved the two-point boundary value problems. Ramos and Rufai [9] implemented a third-derivative, two-step block Falkner approach to solve linear and non-linear BVPs. Some stochastic nonlinear second-order boundary value problems driven by additive noise has been solved by Baccouch [10] using the finite difference method. Numerous physical phenomena are modeled using strongly non-linear BVPs with specific parameter values. In particular, the Bratu's problem is a special type of nonlinear BVP of order two based on eigenvalues. Complex physical and chemical models are frequently described using Bratu's problem in both science and engineering. This problem is employed in a wide range of applications, such as the combustion theory's fuel ignition model, the thermal reaction mechanism, the Chandrasekhar model of universe expansion and chemical reaction theory, radiative heat transmission, and nanotechnology. It has been solved and analysed by different authors using different methods [11]-[16] in the literature.

The Galerkin method is a very well-known method for numerical approximations of various types of problems. Numerous authors come up with various types of problems and solve them numerically using the Galerkin method [17]-[21]. In recent years, the implicit finite difference methodology, also known as the Compact Finite Difference method, has gained popularity. Lele [22] and Mehra & Patel [23] gave a discussion on a variety of ordering schemes in addition to a set of compact finite difference approaches. While Malele *et al.* [24] used the highly precise compact finite difference approach for the solution of BVPs with boundary conditions of Robin type.

In this study, an algorithm is used to improve the accuracy of computation for a class of second-order linear and nonlinear BVPs. For the proposed scheme, we employ the well-known Galerkin and Compact Finite Difference methods. We also introduce the mathematical formulation of 4th-order compact finite difference scheme for the general nonlinear BVPs. We use the Galerkin method to find our first approximation. Then, we use the compact finite difference method to solve the error differential equation in order to improve the accuracy of our first approximation.

## II. GALERKIN METHOD FOR BVPs

The modified Galerkin method is a remarkably potent numerical technique for approximating boundary value problems. The fundamental method employs a trial solution that satisfies all of the problem's boundary conditions.

Let us consider the trial approximation [25] of such a function $f(x)$ of the differential equation (1), given by,

$$\tilde{f}(x) = N_0(x) + \sum_{j=1}^{n} \alpha_j N_j(x) \tag{3}$$

where $\alpha_j$ denotes the undetermined parameters and $N_j(x)$ denotes the basis functions (here Bernstein Polynomials).

The preceding is the generic form of Bernstein polynomials of degree $p$ over the interval $[x_0, x_n]$,

$$B_{p,i}(x) = \binom{p}{i} \frac{(x_n - x)^{p-i}(x - x_0)^i}{(x_n - x_0)^p} \qquad i = 0,1,2,\ldots\ldots p \tag{4}$$

Thus, over the interval [0,1] the Bernstein polynomials of degree $p$ can be represented as,

$$B_{p,i}(x) = \binom{p}{i}(1-x)^{p-i} x^i \qquad i = 0,1,2,\ldots,p \tag{5}$$

These Bernstein polynomials of degree $p$ over the interval [0,1] have some special properties. One of them is, $B_{p,i}(0) = 0 = B_{p,i}(1)$ for $i = 1,2,3,\ldots,p-1$. So, they can be used as the basis of the Galerkin method in the case of Dirichlet boundary conditions.

Now, the Galerkin weighted residual equation becomes [25],

$$\int_a^b \left[ \frac{d^2 \tilde{f}}{dx^2} - g\left(x, \tilde{f}, \frac{d\tilde{f}}{dx}\right) \right] N_i(x) dx = 0 \tag{6}$$

After applying integration by parts to the second derivative part of (6) and then substituting the trial solution (3) into it, we get,

$$\int_a^b \left[ \left( \frac{dN_0}{dx} + \sum_{j=1}^{n} \alpha_j \frac{dN_j}{dx} \right) \frac{dN_i}{dx} + g\left(x, N_0 + \sum_{j=1}^{n} \alpha_j N_j, \frac{dN_0}{dx} + \sum_{j=1}^{n} \alpha_j \frac{dN_j}{dx}\right) N_i \right] dx = \left[ \frac{d\tilde{f}}{dx} N_i \right]_a^b \tag{7}$$

Or equivalently, they can be written as the matrix form as,

$$\sum_{j=1}^{n} \alpha_j K_{ij} = F_i \qquad \text{or} \qquad \mathbf{K}\alpha = \mathbf{F} \tag{8}$$

where,

$$K_{ij} = \int_a^b \left[ \frac{dN_i}{dx} \frac{dN_j}{dx} + Q_1\left(x, N_j, \frac{dN_j}{dx}\right) N_i \right] dx \tag{9}$$

$$F_i = \left[ \frac{d\tilde{f}}{dx} N_i \right]_a^b - \int_a^b \left[ \frac{dN_0}{dx} \frac{dN_i}{dx} + Q_0\left(x, N_0, \frac{dN_0}{dx}\right) N_i \right] dx \tag{10}$$

Solving the system of equations gives us the values of the unknowns, α. In the case of nonlinear BVPs, iteration is required to find the value of these unknown coefficients. Then, the first approximate solution of the BVP is obtained from substituting the values of α's into (3).

## III. COMPACT FINITE DIFFERENCE METHOD FOR BVPS

Firstly, consider the generic form of the differential equation (1) in accordance with the boundary conditions (2). The domain in which the boundary value problem is valid is $[a, b]$.

In order to simplify things further, let us divide the domain $[a, b]$ into $n$ sub-intervals, with the length of each sub-interval being $h$.

It is possible to obtain an implicit numerical approximation of the 1st derivative $f'(x)$ at the mesh points by expressing

$$A f'_{i-1} + f'_i + A f'_{i+1} = \frac{a_1}{2h}(f_{i+1} - f_{i-1}) \tag{11}$$

Here, $A$ and $a_1$ are arbitrarily chosen constants. These constants for 4th order implicit compact finite difference estimations are what we are attempting to discover now. We obtain the following from the Taylor series expansion:

$$f_{i+1} = f_i + h f'_i + \frac{h^2}{2!} f''_i + \frac{h^3}{3!} f'''_i + \frac{h^4}{4!} f^{iv}_i + \frac{h^5}{5!} f^v_i + O(h^6) \tag{12}$$

$$f_{i-1} = f_i - h f'_i + \frac{h^2}{2!} f''_i - \frac{h^3}{3!} f'''_i + \frac{h^4}{4!} f^{iv}_i - \frac{h^5}{5!} f^v_i + O(h^6) \tag{13}$$

Now, subtraction of (13) from (12) becomes,

$$f_{i+1} - f_{i-1} = 2h f'_i + \frac{h^3}{3} f'''_i + \frac{h^5}{60} f^{iv}_i + O(h^7) \tag{14}$$

Again, from the Taylor series expansion:

$$f_{i+1} = f_i + hf_i' + \frac{h^2}{2!}f_i'' + \frac{h^3}{3!}f_i''' + \frac{h^4}{4!}f_i^{iv} + O(h^5) \tag{15}$$

$$f_{i-1} = f_i - hf_i' + \frac{h^2}{2!}f_i'' - \frac{h^3}{3!}f_i''' + \frac{h^4}{4!}f_i^{iv} + O(h^5) \tag{16}$$

Now, adding (16) and (15) we get,

$$f_{i+1} + f_{i-1} = 2f_i + h^2 f_i'' + \frac{h^4}{12}f_i^{iv} + O(h^6) \tag{17}$$

Now,

$$Af_{i-1}' + f_i' + Af_{i+1}' - \frac{a_1}{2h}(f_{i+1} - f_{i-1})$$

$$= A\left[2f_i' + h^2 f_i''' + \frac{h^4}{12}f_i^v + O(h^6)\right] + f_i' - \frac{a_1}{2h}\left[2hf_i' + \frac{h^3}{3}f_i''' + \frac{h^5}{60}f_i^v + O(h^7)\right]$$

$$= (2A + 1 - a_1)f_i' + \left(A - \frac{a_1}{6}\right)h^2 f_i''' + \left(\frac{A}{12} - \frac{a_1}{120}\right)h^4 f_i^v + O(h^6)$$

Setting the coefficient of the first term to zero yields a second-order scheme, whereas setting the coefficients of the first two terms to zero yields a fourth-order scheme.

Setting the first two coefficients equal to zero we get, $A = \frac{1}{4}$, $a_1 = \frac{3}{2}$, then for these values (11) is a two-parameter family of fourth-order scheme with the truncation error,

$$\grave{o}_i = \left(\frac{A}{12} - \frac{a_1}{120}\right) h^4 f_{i+\eta}^v = \frac{1}{120} h^4 f_{i+\eta}^{(5)}$$

So, the first difference equation is,

$$\frac{3}{4h}f_{i-1} - \frac{3}{4h}f_{i+1} + \frac{1}{4}f_{i-1}' + f_i' + \frac{1}{4}f_{i+1}' = 0 \qquad i = 1, n-1 \tag{18}$$

and this difference equation is of convergent order $O(h^4)$, which is equivalent to [26].

For the second difference equation we need our given differential equation. At the interior nodes $x_i$, (1) becomes,

$$f_i'' = g(x, f_i, f_i') \tag{19}$$

In order to solve this problem, we must discover a 4th order approximation of $f_i''$. Now using the Taylor series expansion, we get,

$$f_{i+1} = f_i + hf_i' + \frac{h^2}{2}f_i'' + \frac{h^3}{6}f_i''' + \frac{h^4}{24}f_i^{iv} + \frac{h^5}{120}f_i^v + \frac{h^6}{720}f_i^{vi} + O(h^7) \tag{20}$$

$$f_{i-1} = f_i - hf_i' + \frac{h^2}{2}f_i'' - \frac{h^3}{6}f_i''' + \frac{h^4}{24}f_i^{iv} - \frac{h^5}{120}f_i^v + \frac{h^6}{720}f_i^{vi} + O(h^7) \tag{21}$$

Adding (20) and (21) we get,

$$f_{i+1} + f_{i-1} = 2f_i + h^2 f_i'' + \frac{h^4}{12}f_i^{iv} + \frac{h^6}{360}f_i^{vi} + O(h^8) \tag{22}$$

Again, from the Taylor series expansion we get,

$$f_{i+1}' = f_i' + hf_i'' + \frac{h^2}{2}f_i''' + \frac{h^3}{6}f_i^{iv} + \frac{h^4}{24}f_i^v + \frac{h^5}{120}f_i^{vi} + O(h^6) \tag{23}$$

$$f_{i-1}' = f_i' - hf_i'' + \frac{h^2}{2}f_i''' - \frac{h^3}{6}f_i^{iv} + \frac{h^4}{24}f_i^v - \frac{h^5}{120}f_i^{vi} + O(h^6) \tag{24}$$

Subtracting (24) from (23) we get,

$$f_{i+1}' - f_{i-1}' = 2hf_i'' + \frac{h^3}{3}f_i^{iv} + \frac{h^5}{60}f_i^{vi} + O(h^7) \tag{25}$$

Firstly, multiplying (22) by 4 and (25) by $h$ and then, subtracting (25) form (22) we get,

$$4(f_{i+1} + f_{i-1}) - h(f_{i+1}' - f_{i-1}') = 4\left(2f_i + h^2 f_i'' + \frac{h^4}{12}f_i^{iv}\right.$$

$$+\frac{h^6}{360}f_i^{vi}+O(h^8)\bigg)-h\bigg(2hf_i^{''}+\frac{h^3}{3}f_i^{iv}+\frac{h^5}{60}f_i^{vi}+O(h^7)\bigg) \quad (26)$$

After simplification we get,

$$2h^2 f_i^{''} = 4(f_{i+1}-2f_i+f_{i-1})-h(f_{i+1}^{'}-f_{i-1}^{'})+O(h^6) \quad (27)$$

Or equivalently,

$$f_i^{''} = \frac{2}{h^2}(f_{i+1}-2f_i+f_{i-1})-\frac{1}{2h}(f_{i+1}^{'}-f_{i-1}^{'})+O(h^4) \quad (28)$$

which is an approximation of the second derivative with convergent order $O(h^4)$.

Now, substituting (28) in (19) we get,

$$\frac{2}{h^2}(f_{i+1}-2f_i+f_{i-1})-\frac{1}{2h}(f_{i+1}^{'}-f_{i-1}^{'})=g(x_i,f_i,f_i^{'}) \quad (29)$$

Or equivalently,

$$\frac{2}{h^2}(f_{i+1}-2f_i+f_{i-1})-\frac{1}{2h}(f_{i+1}^{'}-f_{i-1}^{'})-g(x_i,f_i,f_i^{'})=0 \qquad i=1,n-1 \quad (30)$$

This is our 2nd difference equation of convergence order $O(h^4)$ for the nonlinear differential equation (1).

We obtain $(n-1)$ equations for $(n-1)$ interior nodes from both of the difference equations (18) and (30). This indicates that concurrently, they provide a total of $(2n-2)$ equations.

To satisfy the criterion for a unique solution, the number of equations and the number of unknowns must be equal. In addition, we need two additional equations for the Dirichlet and Neumann boundary conditions and four additional equations for the Robin boundary conditions. These extra equations are obtained from the boundary conditions provided.

At node $x=x_0=a$, we get from (19),

$$f_0^{''}=g(x_0,f_0,f_0^{'}) \quad (31)$$

and at node $x=x_n=b$, we get from (19),

$$f_n^{''}=g(x_n,f_n,f_n^{'}) \quad (32)$$

The next step is to establish a fourth-order approximation for $f_0^{''}$ and $f_n^{''}$ in terms of the variables $f_0$, $f_1$, and $f_2$ as well as the variables $f_0^{'}, f_1^{'}$ and $f_2^{'}$. Hence, they can be substituted in (31) and (32) to obtain our desired equations. In order to accomplish this goal, we will need to expand $f_1, f_2, f_1^{'}$, and $f_2^{'}$ in Taylor's series.

Then from Taylor series expression of these terms we obtain,

$$\left.\begin{aligned}f_1 &= f_0+hf_0^{'}+\frac{h^2}{2}f_0^{''}+\frac{h^3}{6}f_0^{'''}+\frac{h^4}{24}f_0^{iv}+\frac{h^5}{120}f_0^{v}+O(h^6)\\ f_2 &= f_0+2hf_0^{'}+2h^2 f_0^{''}+\frac{8h^3}{6}f_0^{'''}+\frac{16h^4}{24}f_0^{iv}+\frac{32h^5}{120}f_0^{v}+O(h^6)\\ f_1^{'} &= f_0^{'}+hf_0^{''}+\frac{h^2}{2}f_0^{'''}+\frac{h^3}{6}f_0^{iv}+\frac{h^4}{24}f_0^{v}+O(h^5)\\ f_2^{'} &= f_0^{'}+2hf_0^{''}+2h^2 f_0^{'''}+\frac{8h^3}{6}f_0^{iv}+\frac{16h^4}{24}f_0^{v}+O(h^5)\end{aligned}\right\} \quad (33)$$

Solving this system for four unknowns $f_0^{''}, f_0^{'''}, f_0^{iv}$, and $f_0^{v}$ we get the value of $f_0^{''}$ as,

$$f_0^{''}=\frac{1}{2h^2}(-23f_0+16f_1+7f_2)-\frac{1}{h}(6f_0^{'}+8f_1^{'}+f_2^{'})+O(h^4) \quad (34)$$

which is the fourth order estimate of the second derivative at the boundary node $x=x_0=a$.

Again, from Taylor series expansion we get,

$$f_{n-1} = f_n - hf_n' + \frac{h^2}{2}f_n'' - \frac{h^3}{6}f_n''' + \frac{h^4}{24}f_n^{iv} - \frac{h^5}{120}f_0^v + O(h^6)$$

$$f_{n-2} = f_n - 2hf_n' + 2h^2 f_n'' - \frac{8h^3}{6}f_n''' + \frac{16h^4}{24}f_n^{iv} - \frac{32h^5}{120}f_n^v + O(h^6)$$

$$f_{n-1}' = f_n' - hf_n'' + \frac{h^2}{2}f_n''' - \frac{h^3}{6}f_n^{iv} + \frac{h^4}{24}f_n^v + O(h^5)$$

$$f_{n-2}' = f_n' - 2hf_n'' + 2h^2 f_n''' - \frac{8h^3}{6}f_n^{iv} + \frac{16h^4}{24}f_n^v + O(h^5)$$

(35)

Solving this system for four unknowns $f_n'', f_n''', f_n^{iv}$, and $f_n^v$ we get the value of $f_n''$ as,

$$f_n'' = \frac{1}{2h^2}(7f_{n-2} + 16f_{n-1} - 23f_n) + \frac{1}{h}(f_{n-2}' + 8f_{n-1}' + 6f_n') + O(h^4) \tag{36}$$

which is the fourth-order estimate of the second derivative at the boundary $x = x_n = b$.

Substituting $f_0''$ in (31) we get,

$$\frac{1}{2h^2}(-23f_0 + 16f_1 + 7f_2) - \frac{1}{h}(6f_0' + 8f_1' + f_2') - g(x_0, f_0, f_0') = 0 \tag{37}$$

This is our first additional equation derived from the left boundary with convergent order $O(h^4)$ for the nonlinear differential equation.

Again, substituting $f_n''$ in (32) we get,

$$\frac{1}{2h^2}(7f_{n-2} + 16f_{n-1} - 23f_n) + \frac{1}{h}(f_{n-2}' + 8f_{n-1}' + 6f_n') - g(x_n, f_n, f_n') = 0 \tag{38}$$

This is our second additional equation derived from the right boundary with convergent order $O(h^4)$ for the nonlinear differential equation.

For Dirichlet and Neumann boundary conditions, the required $2n$ equations come from the two difference equations (18) and (30), and from the two additional equations (37) and (38) at the boundary of the domain. But, for Robin boundary conditions, we required two more additional equations. So, to reduce the number of unknowns, we have to put the boundary conditions into both the difference equations and the equations that come from the boundary. So that, we get a reduced system of $2n$ equations with $2n$ unknowns.

When we evaluate these equations simultaneously, we can find the approximate solutions to any differential equation with different kinds of boundary conditions.

## IV. RESIDUAL CORRECTION

Residual refers to the amount of error remaining after comparing an approximation to the real value of a solution. Boundary value problems can be approximated using a variety of numerical approximation techniques. But there are residuals for each. Residual correction is a technique for minimizing residuals by solving the error differential equation. The residual correction methodology was first presented by Oliveira [27] using the explicit finite difference method. Celik [28] utilizes the Chebyshev series for residual corrections later on. For residual corrections, this section employs the implicit compact finite difference approach.

Consider the differential equation (1) in accordance with its boundary conditions (2) again. As we intend to apply the compact finite difference approach for residual correction, we must first discretize the domain into subdomains. Following this, we are tasked with deriving the differential error equation as well as the boundary conditions that correspond to it from the BVP that was provided.

If $\tilde{f}$ is the approximate solution and $f$ is the exact solution, the residual can be obtained by,

$$\theta = f - \tilde{f} \Rightarrow f = \tilde{f} + \theta \tag{39}$$

Substituting this result into (1) we get,

$$\frac{d^2}{dx^2}(\tilde{f} + \theta) = g\left(x, \tilde{f} + \theta, \frac{d\tilde{f}}{dx} + \frac{d\theta}{dx}\right) \tag{40}$$

Or equivalently, this can be written as,

$$\frac{d^2\theta}{dx^2} - M_1\left(x, \theta, \frac{d\theta}{dx}\right) = -\frac{d^2\tilde{f}}{dx^2} + M_2\left(x, \tilde{f}, \frac{d\tilde{f}}{dx}\right) \tag{41}$$

which is the differential equation of the error which follows the subsequent boundary conditions,

$$\left.\begin{array}{l} a_1\tilde{f}(a) + a_1\theta(a) + b_1\left.\frac{d\tilde{f}}{dx}\right|_{x=a} + b_1\left.\frac{d\theta}{dx}\right|_{x=a} = \lambda_1 \\[1em] a_2\tilde{f}(b) + a_2\theta(b) + b_2\left.\frac{d\tilde{f}}{dx}\right|_{x=b} + b_2\left.\frac{d\theta}{dx}\right|_{x=b} = \lambda_2 \end{array}\right\} \tag{42}$$

But the approximate solution satisfies the given boundary conditions. So, the boundary conditions for the error differential equations becomes,

$$\left.\begin{array}{l} a_1\theta(a) + b_1\left.\frac{d\theta}{dx}\right|_{x=a} = 0 \\[1em] a_2\theta(b) + b_2\left.\frac{d\theta}{dx}\right|_{x=b} = 0 \end{array}\right\} \tag{43}$$

We solve the BVP (41) in accordance with the boundary conditions (43) by 4$^{th}$ order compact finite difference method. Then, the updated approximate value becomes,

$$(\text{Updated Approximation}) = (\text{Weighted Residual Value}) + (\text{Error Value}) \tag{44}$$

## V. CONVERGENCE ANALYSIS

In this section we discuss the convergence of our proposed scheme for solving the linear and nonlinear BVPs (1)-(2).

Let $\Omega = C^l([x_0, x_n] \in R)$ be the linear space of real-valued functions that are $l$ times differentiable on $D = [x_0, x_n]$. Suppose that,

$$<\zeta_1, \zeta_2> = \int_D w_0(x)\zeta_1(x)\zeta_2(x)dx \tag{45}$$

be the $L^2$ inner product on $\Omega$ for some sufficiently smooth weight function $w_0$ that induces the $L^2$ norm,

$$|\zeta_1|^2 = \int_D w_0(x)\zeta_1^2(x)dx \tag{46}$$

for which $\Omega$ is an infinite dimensional Hilbert space. Assume $V = \{B\_i | i = 1,2,3, \dots \dots\}$ be a Schauder basis of $\Omega$, where $B_i$'s are Bernstein polynomials on $D = [x_0, x_n]$. Let us begin with an approximation subspace $\Omega^N$ spanned by $\{\psi_1, \psi_2, \psi_3 \dots \dots \dots, \psi_N\}$ that satisfies appropriate boundary conditions. The Galerkin weighted residual equation is given by $< R\left(x, \tilde{f}(x)\right), \psi_j >= 0$ where $\tilde{f}(x) = N_0(x) + \sum_{j=1}^n \alpha_j N_j(x)$ is a trial solution of the differential equation. In particular, the residual function $R\left(x, \tilde{f}(x)\right)$ is orthogonal to each function of the basis of the approximate subspace $\Omega^N$. It is well known that Bernstein polynomials can approximate any continuous function with arbitrary precision. If the dimension of the subspace $\Omega^N$ is infinitely large, the residual function $R\left(x, \tilde{f}(x)\right)$ is orthogonal to each Bernstein polynomial which immediately implies that the residual $R\left(x, \tilde{f}(x)\right)$ is orthogonal to any continuous function in $\Omega$. A function that is orthogonal to any other functions in the space is necessarily the zero function.

For the error differential equation in accordance with the error boundary conditions, we use 4th order compact finite difference method.

**Definition 1**: [29] Let us consider $\Psi_h$ as the solution of the discretized equation $\mathcal{L}_h\Psi_h = b_h$ which converges to the solution $\psi$ of the given differential equation $\mathcal{L}\psi = b$ if $||\psi_h - \Psi_h|| \to 0$ as $h \to 0$.

Moreover, if for a positive constant $k$ and another positive constant $M_0 > 0$ which is independent of $k$ in the sense that, $\|\psi_h - \Psi_h\| \leq M_0 h^k$ in this case, it is stated that the discretized equation has $k$-th order precision with convergence order $h^k$.

**Definition 2:** [29] The discretized Equation is stable if there exists $h_0 > 0$ and $\delta > 0$ in such a way that for any $h < h_0$ and any $\epsilon_h \in \mathcal{B}_h$ satisfying $\|\epsilon\| <$, the perturbed difference equation $\mathcal{L}_h w_h = b_h + \delta \epsilon_h$ has unique solution $w_h$, satisfying $\|w_h - \Psi_h\| \leq M \|\epsilon_h\|$ where $\Psi_h$ is the solution of the unperturbed difference equation and $M > 0$ is independent of $h$.

**Theorem 1:** [29] *If the discretized equation is stable and also consistent (with order $h^k$) with the given differential equation, then, the solution of the discretized equation $\Psi_h$ converges to the solution $\psi$ and which satisfies $\|\psi_h - \Psi_h\| \leq M M_1 h^k$, where $M$ and $M_1$ are certain constants. Alternatively, the order in which the difference scheme approaches the continuous problem corresponds to the accuracy order of the difference scheme.*

**Proof:** [29] Since the difference scheme is consistent, it is established that $\|\delta b_h\| \leq M_1 h^k \to 0$ as $h \to 0$. This means that a grid $\mathcal{D}_h$ can be constructed in such way that $h < h_0$ and $\delta b_h <$ as in Definition 2, and

$$L_h \Psi_h = b_h + \delta b_h \tag{47}$$

Therefore, $\Psi_h$ satisfies the conditions for stability. Due to this,

$$\|\psi_h - \Psi_h\| \leq M \|\delta b_h\| \leq M(M_1 h^k) \tag{48}$$

Hence, the discrete solution $\Psi_h$ converges to the continuous solution $\psi$ of order $h^k$.

## VI. RESULTS AND DISCUSSIONS

In this section, we will apply our proposed scheme to some second-order linear and nonlinear problems and calculate the accuracy and stability of our proposed scheme. We compute the $L_\infty$ norm by comparing the approximate solutions with the exact solutions as follows,

$$L_\infty = \max_{1 \leq x \leq n} |\tilde{f}_i - f_i| \tag{49}$$

We also compute the Convergence Rate ($\mathcal{CR}$) as follows [24],

$$\mathcal{CR} = \frac{\log\left(\frac{\epsilon_1}{\epsilon_2}\right)}{\log\left(\frac{h_1}{h_2}\right)} \tag{50}$$

where, $\epsilon_1$ and $\epsilon_2$ are the $L_\infty$ for grid size $h_1$ and $h_2$, respectively.

**Problem 1**

Let us consider the linear DE in conjunction with Neumann boundary conditions [6],[30],

$$\frac{d^2 f}{dx^2} + f = -1 \quad \text{with} \quad f'(0) = \frac{1 - \cos(1)}{\sin(1)}, f'(1) = -\frac{1 - \cos(1)}{\sin(1)} \tag{51}$$

with the exact solution, $\qquad f(x) = \cos(x) + \frac{1 - \cos(1)}{\sin(1)} \sin(x) - 1 \tag{52}$

In the first stage, the BVP is solved numerically using the modified Galerkin approach with the aid of the trial solution as (3). For residual correction we need the error differential equation with error BCs. The error BVP for the given differential equation is,

$$\frac{d^2 \theta}{dx^2} + \theta = -\frac{d^2 \tilde{f}}{dx^2} - \tilde{f} - 1 \quad \text{with} \quad \theta'(0) = 0, \theta'(1) = 0 \tag{53}$$

The following table shows $L_\infty$ and $\mathcal{CR}$ generated using the proposed method for different grid sizes $h$ for residual corrections.

TABLE I: $L_\infty$ AND CONVERGENCE RATE ($\mathcal{CR}$) FOR PROBLEM 1

|  | h | $L_\infty$ | $\mathcal{CR}$ | Reference Results ($L_\infty$) |
|---|---|---|---|---|
|  | 0.1 | $2.0797 \times 10^{-08}$ |  | [6], 2011: |
|  | 0.05 | $1.5311 \times 10^{-09}$ | 3.7637 | Bernoulli Pol. (5th degree) |
| Bernstein | 0.025 | $1.3140 \times 10^{-10}$ | 3.5425 | $1.173569 \times 10^{-06}$ |
| Polynomials | 0.0125 | $9.3694 \times 10^{-12}$ | 3.8099 | Bernoulli Pol. (7th degree) |
| of degree 4 | 0.01 | $3.9338 \times 10^{-12}$ | 3.8892 | $1.2808 \times 10^{-09}$ |
|  | 0.005 | $2.5785 \times 10^{-13}$ | 3.9313 | Bernoulli Pol. (10th degree) |
|  | 0.0025 | $1.6570 \times 10^{-14}$ | 3.9599 | $1.0131 \times 10^{-14}$ |

From Table I we can observe that we obtain high accuracy by using Bernstein polynomials of degree 4 only in our prescribed scheme. Whereas, in [6], Islam & Shirin implemented the modified Galerkin approach and the $L_\infty$ obtained using Bernoulli polynomials of degree 10 is $10^{-14}$. So, we can say that, in our proposed scheme, we can attain better accuracy by using fewer polynomials and then correcting the residuals of the approximation.

**Problem 2**

Let us consider the following nonlinear DE in conjunction with Dirichlet boundary conditions [6],

$$\frac{d^2 f}{dx^2} + \frac{1}{8} f \frac{df}{dx} = \left(4 + \frac{1}{4} x^3\right) \quad \text{with} \quad f(1) = 17, f(3) = \frac{43}{3} \tag{54}$$

with the exact solution, 
$$f(x) = x^2 + \frac{16}{x} \tag{55}$$

Equivalent BVP on [0,1] becomes,

$$\frac{d^2 f}{dx^2} + \frac{1}{4} f \frac{df}{dx} = 16 + (2x+1)^3 \quad \text{with} \quad f(0) = 17, f(1) = \frac{43}{3} \tag{56}$$

In the first stage, the BVP is solved numerically using the modified Galerkin approach with the aid of the trial solution as (3). For residual correction we need the error differential equation with error BCs. The error BVP for the given differential equation is,

$$\frac{d^2 \theta}{dx^2} + \frac{1}{4} \tilde{f} \frac{d\theta}{dx} + \frac{1}{4} \theta \frac{d\tilde{f}}{dx} + \frac{1}{4} \theta \frac{d\theta}{dx} = 16 + (2x+1)^3 - \frac{d^2 \tilde{f}}{dx^2} - \frac{1}{4} \tilde{f} \frac{d\tilde{f}}{dx} \quad \text{with} \quad \theta(0) = \theta(1) = 0 \tag{57}$$

The following table shows $L_\infty$ and $\mathcal{CR}$ generated using the proposed method for different grid sizes $h$ for residual corrections.

TABLE II: $L_\infty$ AND CONVERGENCE RATE ($\mathcal{CR}$) FOR PROBLEM 2

|  | h | $L_\infty$ | $\mathcal{CR}$ | Reference Results ($L_\infty$) |
|---|---|---|---|---|
|  | 0.1 | $5.1458 \times 10^{-04}$ |  | [6], 2011: |
|  | 0.05 | $5.1111 \times 10^{-05}$ | 3.3315 | Bernoulli Pol. (10th degree) |
| Bernstein | 0.025 | $4.2499 \times 10^{-06}$ | 3.5883 | $1.4449 \times 10^{-05}$ |
| Polynomials | 0.0125 | $3.1017 \times 10^{-07}$ | 3.7763 | Bernoulli Pol. (15th degree) |
| of degree 4 | 0.01 | $7.6237 \times 10^{-08}$ | 3.8548 | $6.3806 \times 10^{-08}$ |
|  | 0.005 | $8.7611 \times 10^{-09}$ | 3.9048 |  |
|  | 0.0025 | $5.6623 \times 10^{-10}$ | 3.9516 |  |

From Table II we can observe that we obtain high accuracy by using Bernstein polynomials of degree 4 only in our prescribed scheme. Whereas, in [6], Islam & Shirin implemented the modified Galerkin approach and the $L_\infty$ obtained using Bernoulli polynomials of degree 15 is $10^{-08}$. So, we can say that, in our proposed scheme, we can attain better accuracy by using fewer polynomials and then correcting the residuals of the approximation.

**Problem 3**

Let us consider the following nonlinear DE in conjunction with Robin boundary conditions [6],[31],

$$\frac{d^2 f}{dx^2} = \frac{1}{2}(1 + x + f)^3 \quad \text{with} \quad f'(0) - f(0) = -1/2, \, f'(1) + f(1) = 1 \tag{58}$$

with the exact solution, 
$$f(x) = \frac{2}{2-x} - x - 1 \tag{59}$$

In the first stage, the BVP is solved numerically using the modified Galerkin approach with the aid of the trial solution as (3). For residual correction we need the error differential equation with error BCs. The error BVP for the given differential equation is,

$$\frac{d^2\theta}{dx^2} - \left(\frac{3}{2}(1+x)^2 + 3(1+x)\tilde{f} + \frac{3}{2}\tilde{f}^2\right)\theta - \left(\frac{3}{2}(1+x) + \frac{3}{2}\tilde{f}\right)\theta^2 - \frac{1}{2}\theta^3 = -\frac{d^2\tilde{f}}{dx^2} + \frac{1}{2}(1+x+\tilde{f})^3 \quad (60)$$

with the boundary conditions, $\theta'(0) - \theta(0) = -1/2$, $\theta'(1) + \theta(1) = 1$.

The following table shows $L_\infty$ and $\mathcal{CR}$ generated using the proposed method for different grid sizes $h$ for residual corrections.

TABLE III: $L_\infty$ AND CONVERGENCE RATE ($\mathcal{CR}$) FOR PROBLEM 3

| | h | $L_\infty$ | $\mathcal{CR}$ | Reference Results ($L_\infty$) | |
|---|---|---|---|---|---|
| | | | | [6], 2011: | [31], 2022 |
| | 0.1 | $7.1799 \times 10^{-06}$ | | Bernoulli Pol. ($8^{th}$ degree) | Legendre Pol. ($8^{th}$ degree) |
| | 0.05 | $2.8731 \times 10^{-07}$ | 4.6433 | $5.52 \times 10^{-07}$ | $2.65 \times 10^{-09}$ |
| Bernstein | 0.025 | $1.1124 \times 10^{-08}$ | 4.6909 | Bernoulli Pol. ($10^{th}$ degree) | Legendre Pol. ($10^{th}$ degree) |
| Polynomials of | 0.0125 | $8.9334 \times 10^{-10}$ | 3.6383 | $6.99 \times 10^{-09}$ | $7.77 \times 10^{-11}$ |
| degree 4 | 0.01 | $3.8393 \times 10^{-10}$ | 3.7846 | | Legendre Pol. ($11^{th}$ degree) |
| | 0.005 | $2.6363 \times 10^{-11}$ | 3.8643 | | $8.15 \times 10^{-12}$ |
| | 0.0025 | $1.7256 \times 10^{-12}$ | 3.9333 | | |

From Table III we can observe that we obtain high accuracy by using Bernstein polynomials of degree 4 in our proposed scheme. Whereas, in [6], Islam & Shirin implemented the modified Galerkin approach using Bernoulli polynomials of degree 10 and the obtained $L_\infty$ is $10^{-9}$. In [31], Sohel et al. implemented the residual correction approach and the $L_\infty$ using Legendre polynomials of degree 11 is $10^{-12}$. So, we can say that, we can attain better accuracy by using fewer polynomials and then correcting the residuals of the approximation using our present approach.

**Problem 4**

Let us consider the following strongly nonlinear Bratu's problem [12],[15],[32]-[35],

$$\frac{d^2 f}{dx^2} + \lambda e^{f(x)} = 0, \text{ with } f(0) = 0, f(1) = 0 \quad (61)$$

The term strongly non-linear is employed because the non-linearity is generated by exponential term.

This BVP has an exact solution, $$f(x) = -2\ln\left[\frac{\cosh\left(x - \frac{1}{2}\right)\frac{\beta}{2}}{\cosh\left(\frac{\beta}{4}\right)}\right] \quad (62)$$

where, $\beta$ satisfies, $\sqrt{2\lambda}\cosh\left(\frac{\beta}{4}\right) = \beta$. There are zero, one, and two solutions to the Bratu's problem when $\lambda > \lambda_c$, $\lambda = \lambda_c$ and $\lambda < \lambda_c$, respectively, where $\lambda_c = 3.513830719$.

In the first stage, the BVP is solved numerically using the modified Galerkin approach with the aid of the trial solution as (3). For residual correction we need the error differential equation with error BCs. The error BVP for the given differential equation is,

$$\frac{d^2\theta}{dx^2} + \lambda e^{\tilde{f}} e^\theta = -\frac{d^2\tilde{f}}{dx^2} \quad \text{with} \quad \theta(0) = 0, \theta(1) = 0 \quad (63)$$

The following table shows $L_\infty$ and $\mathcal{CR}$ generated using the proposed method for different grid sizes $h$ for residual corrections. Firstly, representing the results for $\lambda = 1$.

TABLE IV: $L_\infty$ AND CONVERGENCE RATE ($\mathcal{CR}$) FOR PROBLEM 4 FOR $\lambda = 1$

| | h | $L_\infty$ | $\mathcal{CR}$ | Reference Results ($L_\infty$) |
|---|---|---|---|---|
| | 0.1 | $8.2287 \times 10^{-08}$ | | [12], 2021: $6.187 \times 10^{-05}$ |
| Bernstein | 0.05 | $7.6868 \times 10^{-09}$ | 3.4202 | [15], 2004: $1.348 \times 10^{-05}$ |
| Polynomials | 0.025 | $6.2340 \times 10^{-10}$ | 3.6242 | [32], 2011: $9.048 \times 10^{-07}$ |
| of | 0.01 | $1.8176 \times 10^{-11}$ | 3.8580 | [34], 2019: $9.90 \times 10^{-08}$ |
| degree 4 | 0.005 | $1.1832 \times 10^{-12}$ | 3.9413 | [35], 2019: $6.41 \times 10^{-13}$ |
| | 0.0025 | $7.6439 \times 10^{-14}$ | 3.9522 | |

From Table IV we can observe that we obtain the accuracy of $10^{-14}$ by using Bernstein polynomials of degree 4 in our present scheme for $\lambda = 1$. Whereas, for $\lambda = 1$, in [12], Mustafa et al. implemented Subdivision collocation method and in [15], Khuri implemented Laplace transformation method and the obtained $L_\infty$ in both of the cases is $10^{-05}$, in [32], Abbasbandy et al. implemented the Lie Group Shooting method and the obtained $L_\infty$ is of $10^{-0}$, in [34] Ala O et al. implemented the Quantic B-spline approach and the obtained $L_\infty$ is of $10^{-08}$ and in [35] Roul & Thula implemented the B-spline Collocation method and the obtained $L_\infty$ is of $10^{-13}$. Then, representing the results for $\lambda = 2$,

TABLE V: $L_\infty$ AND CONVERGENCE RATE ($\mathcal{CR}$) FOR PROBLEM 4 FOR $\lambda = 2$

|  | $h$ | $L_\infty$ | $\mathcal{CR}$ | Reference Results ($L_\infty$) |
|---|---|---|---|---|
|  | 0.1 | $9.5936 \times 10^{-07}$ |  | [33], 2000: $1.52 \times 10^{-02}$ |
| Bernstein | 0.05 | $1.1273 \times 10^{-07}$ | 3.0892 | [12], 2021: $4.939 \times 10^{-04}$ |
| Polynomials | 0.025 | $8.6011 \times 10^{-09}$ | 3.7122 | [32], 2011: $5.709 \times 10^{-06}$ |
| of | 0.01 | $2.4317 \times 10^{-10}$ | 3.8616 | [34], 2019: $1.60 \times 10^{-06}$ |
| degree 4 | 0.005 | $1.5668 \times 10^{-11}$ | 3.9561 | [36], 2015: $1.98 \times 10^{-09}$ |
|  | 0.0025 | $9.9537 \times 10^{-13}$ | 3.9764 |  |

From Table V we can observe that we obtain the accuracy of $10^{-13}$ by using Bernstein polynomials of degree 4 in our present scheme for $\lambda = 2$. Whereas, for $\lambda = 2$, in [12], Mustafa et al. implemented the subdivision collocation method and the obtained $L_\infty$ is of $10^{-04}$, in [32], Abbasbandy et al. implemented the Lie Group Shooting method and the obtained $L_\infty$ is of $10^{-06}$, in [33], Deeba et al. implemented the Decomposition algorithm and the obtained $L_\infty$ is of $10^{-02}$, in [34], Ala O et al. implemented the Quantic B-spline method and the obtained $L_\infty$ is of $10^{-06}$ and in [36], Farzana & Islam implemented the Chebyshev-Legendre collocation method and the obtained $L_\infty$ is of $10^{-0}$. So, we can say that, in our proposed scheme, we can attain better accuracy by using fewer polynomials and then correcting the residuals.

## VII. Conclusion

In this research, we have obtained numerical solutions for both linear and nonlinear BVPs, followed by residual corrections. Using the weighted residual technique, a numerical solution has been generated in the first part. Then, using the implicit compact finite difference technique, updated approximations were established. We have developed the formulation of the compact finite difference scheme for the nonlinear differential equations both at the interior node and the boundary with Dirichlet, Neumann, and Robin boundary conditions and established the convergence and stability of our numerical solutions. Our research has shown the effectiveness of the residual correction technique in improving the accuracy of the approximations. We have compared our results with those published in the literature and demonstrated the superiority of our approximations in terms of accuracy. Overall, this research offers valuable insights into the numerical solution of boundary value problems and provides a practical method for achieving high-precision results for various applications. The proposed method can be applied to higher-order boundary value problems as well as to partial differential equations.